\documentclass{article}
\usepackage[utf8]{inputenc}
\usepackage[english]{babel}
\usepackage{amssymb}
\usepackage{amsmath}
\usepackage{csquotes}
\usepackage{cite}

\usepackage{hyperref}
\usepackage{tikz-cd}
\usepackage{tikz}
\usepackage{mathtools} 
\usepackage{extarrows} 
\usetikzlibrary{cd}
\hypersetup{
    colorlinks=true,
    linkcolor=blue,
    filecolor=magenta,      
    urlcolor=cyan,
}

\usepackage{amsthm}
\usepackage{algorithm}
\usepackage[noend]{algpseudocode}
\usepackage{enumitem}
\usepackage{authblk}
\newtheorem{theorem}{Theorem}[section]
\newtheorem{lemma}{Lemma}[section]
\newtheorem{proposition}{Proposition}[section]
\newtheorem{corollary}{Corollary}[section]
\theoremstyle{definition}

\newtheorem{remark}{Remark}[section]

\title{On possible sums from multiset of mutually divisible natural numbers}
\author[1]{Yizhou Guo \thanks{guoyiz@yandex.ru} \thanks{st095712@student.spbu.ru}}
\affil[1]{Saint-Petersburg State University, Department of Mathematics and Computer Science}

\begin{document}
\maketitle
\begin{abstract}
Let $A$ be a finite multiset of $\mathbb{N}$ such for any two $a_i,a_j \in A$, $a_i | a_j$ or $a_j | a_i$. We describe the structure of the set $\mathrm{span}(A)$ of all sums of submultisets of $A$, and in particular give a criterion of $\mathrm{span}(A)=\mathrm{span}(B)$ for two multisets $A,B$.
\end{abstract}
\section{Introduction}
Let $A$ be a finite multiset wherein every element of $A$ is a positive integer. Denote by $\mathrm{sum}(A)$ the sum of all elements of $A$ (respecting multiplicity), and by $\mathrm{span(A)}$ the set (not a multiset) of all sums $\mathrm{sum}(C)$, where $C$ is a submultiset
of $A$.  

This setting describes a finite collection of tokens of positive integer costs, and $\mathrm{span(A)}$ means the set of values which may be obtained using this collection of tokens.

In the general case describing the structure of $\mathrm{span(A)}$ seems difficult. Hereafter we restrict to the case where the numbers which occur in this multiset are of the form $d_0, d_0d_1, \ldots, d_0d_1\cdots d_k$, wherein without loss of generality $d_0=1$ and $d_i>1$ for $i>0$. From now on, we shall denote $p_i = d_0d_1\ldots d_i$.

A more convenient notation would be $A:=(a_0,a_1,\ldots,a_k)$, where $a_i \in \mathbb{Z}_{\ge 0}$ and $a_k \neq 0$, where $a_i$ is the multiplicity of $p_i$ in $A$. We will from here on call this a \textbf{$\{d_i\}$-collection} and let $\mathcal{A}_{\{d_i\}}$ denote the family of $\{d_i\}$-collections. Let $\mathrm{comb}(A)$ denote $\{(c_0,c_1,\ldots,c_k) \in  \mathbb{Z}^{k+1} : 0\leq c_i \leq a_i \}$.

Then $\mathrm{span(A)} = \{\sum_{i=0}^k c_ip_i : (c_0,c_1,\ldots,c_k) \in \mathrm{comb}(A) \}$. 

Here, we will describe and prove the correctness of an algorithm to decide for arbitrary $A,B \in \mathcal{A}_{\{d_i\}}$ whether or not $\mathrm{span(A)}=\mathrm{span(B)}$, as well as an algorithm to decide for arbitrary $n \in \mathbb{Z}_{\ge 0}$ if $n \in \mathrm{span}(A)$. 

An \textbf{elementary exchange} of $(a_0,a_1,\ldots,a_k)$ picks some $i$ such that $a_i \geq d_i$, sets $a_i \leftarrow a_i-d_i$ and $a_{i+1} \leftarrow a_{i+1}+1$. In other words, $d_i$ tokens of cost $p_i$ are replaced with one token of cost $p_{i+1}$. We let $e_i$ denote the elementary exchange wherein the $i$-th place is exchanged. 

For $A=(a_0,a_1,\ldots,a_k),B=(b_0,b_1,\ldots,b_k) \in \mathcal{A}_{\{d_i\}}$, we say that $A \preceq B$ if and only if $a_i \leq b_i$ for all $i$.

\begin{remark}
One can consider the special case where only $1$ and powers of $10$ are allowed. From this, one can generalize to powers for $d$ for arbitrary $d\geq 2$. We can call the latter a \textbf{$d$-collection} and use $\mathcal{A}_d$ to denote the family of $d$-collections. For simplicity, one can first solve this problem for $d$-collections and then generalize to $\{d_i\}$-collections.
\end{remark}

\begin{proposition} \label{prop:span_preserves_partial_order}
    If $A \preceq B$, then $\mathrm{span}(A) \subseteq \mathrm{span}(B)$.
\end{proposition}
\begin{proof}
    That $A \preceq B$ implies $\mathrm{comb}(A) \subseteq \mathrm{comb}(B)$, from which the result follows.
\end{proof}
\begin{proposition} \label{prop:span_of_exchange_subset_of_span}
    \[\mathrm{span}(e_i(A)) \subseteq \mathrm{span}(A)\]
\end{proposition}
\begin{proof}
Instead of using the new token of cost $p_{i+1}$ you could use $d_i$ tokens of cost $p_i$ which it was changed for.
\end{proof}
\section{Condition for span invariance with respect to elementary exchanges}
\begin{proposition} 
\label{prop:invariance_under_elementary_exchange}
    For $A=(a_0,a_1,\ldots,a_k)$, if $a_i > 2(d_i-1)$, then $\mathrm{span}(A) = \mathrm{span}(e_i(A))$.
\end{proposition}
\begin{proof} It suffices to check that  $\mathrm{span}(A)\subset \mathrm{span}(e_i(A))$, since the opposite inclusion is proved
in Proposition \ref{prop:span_of_exchange_subset_of_span}.
It suffices to consider an arbitrary submultiset $B\subset A$ and prove that $\mathrm{sum}(B)\in \mathrm{span}(e_i(A))$. If $B$ contains at most $d_i-1$ tokens $p_i$, it is also a submultiset of $e_i(A)$, since $e_i(A)$ contains more than $2(d_i-1)-d_i=d_i-2$ tokens $p_i$. If $B$ contains at least $d_i$ tokens $p_i$, replace $d_i$ such tokens by a token $p_{i+1}$ and get a submultiset $e_i(B)\subset e_i(A)$ with $\mathrm{sum}(e_i(B))=\mathrm{sum}(B)$. Thus in both cases $\mathrm{sum}(B)\in \mathrm{span}(e_i(A))$, as needed.
\end{proof}
From now on, we will call an elementary exchange at an $i$ where the hypothesis of the above proposition is satisfied to be a \textbf{proper elementary exchange}.
\begin{algorithm}
\caption{Proper elementary exchange sequence}\label{alg:proper_elem_exchange_seq}
\begin{algorithmic}
\Procedure{RunProperElementaryExchanges}{$A$}
\Require for all $i, A[i]=a_i\cdot p_i, 0 \leq a_i < d_i$
\Ensure for all $i, A[i]\leq 2(d_i-1)$
\State $i \gets 0$
\While{exists untraversed non-zero element of $A$}
\While{$A[i] > 2(d_i-1)$}
\State $A[i] \gets A[i] - d_i$
\State $A[i+1] \gets A[i+1]+1$
\EndWhile
\State $i \gets i+1$
\EndWhile
\Return A
\EndProcedure
\end{algorithmic}
\end{algorithm}
\begin{corollary} \label{cor:span_of_collection_equals_span_of_its_normalization}
    Let $B$ be the $\{d_i\}$-collection returned by running Algorithm \ref{alg:proper_elem_exchange_seq} on arbitrary $\{d_i\}$-collection $A$. Then, $\mathrm{span}(A)=\mathrm{span}(B)$.
\end{corollary}
\begin{proof}
    Let $A,A_1,A_2,\ldots,A_m,B$ be the sequence of $\{d_i\}$-collections obtained from running the algorithm on $A$. This then follows from Proposition \ref{prop:invariance_under_elementary_exchange} applied to adjacent elements of this sequence along with transitivity of set equality.
\end{proof}

\begin{proposition}
The converse of Proposition \ref{prop:invariance_under_elementary_exchange} is false.    
\end{proposition}
\begin{proof}
A simple counterexample would be for the $2$-collection $A = (3,2,0)$, with respect to which there is an elementary exchange that results in $A'=(3,0,1)$ in which case $2(d-1)=2 \geq d-1$ but all $7$ values can be assumed by both $A$ and $A'$.
\end{proof}

We say that $\{d_i\}$-collection $A$ is \textbf{normal} if and only if $a_i \leq 2(d-1)$ for all $i$ and that $\{d_i\}$-collection $A$ is \textbf{$j$-normal} if and only if $d_i-1 \leq a_i \leq 2(d_i-1)$ for all $i<j$ and $a_j < d_i-1$. We also define the \textbf{normalization} and \textbf{$j$-normalization} of a $\{d_i\}$-collection $A$ in the same manner per Algorithm \ref{alg:proper_elem_exchange_seq} (with the latter running a modified version that does not exchange at any place $\geq j$), which we denote with $\mathrm{norm}(A)$ and $\mathrm{norm}_j(A)$ respectively. $[\mathrm{norm}(A)]_i$ will denote the value at its $i$th place (such subscripting can be applied to any expression evaluating to a $\{d_i\}$-collection type).
\begin{proposition} \label{prop:same_normalization_implies_same_span}
    If $\mathrm{norm}(A)=\mathrm{norm}(B)$, then $\mathrm{span}(A)=\mathrm{span}(B)$.
\end{proposition}
\begin{proof}
    Follows from Proposition \ref{cor:span_of_collection_equals_span_of_its_normalization} along with transitivity of set equivalence.
\end{proof}
\begin{lemma} \label{lemma:condition_for_assuming_all_values}
    If a $\{d_i\}$-collection $A=(a_0,a_1,\ldots,a_k)$ is such that $a_i \geq d_i-1$ for all $i<k$, then $\mathrm{span}(A)=
    \{0,1,\ldots,\mathrm{sum}(A)\}$.
    \end{lemma}
\begin{proof}
    We prove this by induction by $\mathrm{sum}(A)$. The base case follows trivially. The inductive hypothesis applied to $A\setminus \{p_k\}$ implies that all elements  up to $M:=\mathrm{sum}(A)-p_k$ belong to $\mathrm{span}(A)$. Thus so are all elements from $\mathrm{sum}(A)-M$ to $\mathrm{sum}(A)$, since the set $\mathrm{span}(A)$ enjoys a symmetry $x\to \mathrm{sum}(A)-x$. It remains to note that $\mathrm{sum}(A)\geq p_k+\sum_{i=0}^{k-1} (d_i-1)p_i=2p_k-1$,
    thus $\mathrm{sum}(A)-M=p_k\leq \mathrm{sum}(A)-p_k+1=M+1$.
\end{proof}
\begin{proposition} \label{prop:smallest_value_not_in_span}
    If a $\{d_i\}$-collection $A=(a_0,a_1,\ldots,a_k)$ is $j$-normal, then the smallest value not in $\mathrm{span}(A)$ is 
    $X:=1+\sum_{i=0}^{j} a_ip_i$, which is necessarily $< p_{j+1}$.
\end{proposition}
\begin{proof}
    That $A$ is $j$-normal implies that
    \begin{eqnarray*}
    X-1=a_jp_j+\sum_{i=0}^{j-1} a_ip_i &\leq& 2(p_j-1)+a_jp_j\\
    &\leq& 2(p_j-1)+(d_j-2)p_j = p_{j+1}-2.
    \end{eqnarray*}
    By Lemma \ref{lemma:condition_for_assuming_all_values}, any $n \leq \sum_{i=0}^j a_ip^i=X-1$ is in $\mathrm{span}(A)$. It remains to prove that $X\notin \mathrm{span}(A)$. Since $X<p_{j+1}$, we can use only tokens not exceeding $p_j$ for collecting $X$, but the sum of such tokens in $A$ equals $X-1$, thus indeed $X\notin \mathrm{span}(A)$.
\end{proof}

\section{Decomposition of normal $\{d_i\}$-collections}
We now describe an algorithm for decomposing an arbitrary normal $\{d_i\}$-collection $A$ to submultisets which partition $A$. Call an index $j\leq k$ critical if either $j=k$ or $a_j<d_i-1$. Let $j_1<j_2<\ldots<j_s=k$ be all critical indices. Define the $\{d_i\}$-collections $A_1,\ldots,A_s$ as follows: $A_i$ consists of those tokens $p_t$ from $A$, for which $j_{i-1}<t\leq j_i$ (for $i=1$, we put $j_0=-1$ by agreement). If $s=1$, we say that $A$ is an \emph{irreducible} $\{d_i\}$-collection.
\begin{theorem}
In the above notation, 
$$
\mathrm{span}(A)=\oplus_{i=1}^s p_{j_{i-1}+1}\cdot \{0,1,\ldots, p_{-1-j_{i-1}}\mathrm{sum}(A_j)\},
$$
where the $\oplus$ sign means the direct sum: any element in 
$\mathrm{span}(A)$ is uniquely represented as $x_1+\ldots+x_s$ for 
$x_i\in p_{j_{i-1}+1}\cdot \{0,1,\ldots, p_{-1-j_{i-1}}\mathrm{sum}(A_j)\}$.
\end{theorem}
\begin{proof}
It follows from Lemma \ref{lemma:condition_for_assuming_all_values} and scaling that
$$\mathrm{span}(A_i)=p_{j_{i-1}+1}\cdot \{0,1,\ldots, p_{-1-j_{i-1}}\mathrm{sum}(A_j)\}.$$
Also, by Proposition \ref{prop:smallest_value_not_in_span} we get $\mathrm{sum}(A_i)\leq p_{j_i+1}-p_{j_{i-1}+1}$ for all $i<s$. Since any element $x\in \mathrm{span}(A)$ may be represented as $x=\sum x_i$ for $x_i\in \mathrm{span}(A_i)$, it remains to prove the direct sum condition. Assume on the contrary that $x=\sum x_i=\sum y_i$ for $x_i,y_i\in \mathrm{span}(A_i)$ but there exists $i$ for which $x_i\ne y_i$. Choose minimal such $i$, obviously we have $i<s$. Then $x_i-y_i=(y_{i+1}+\ldots+y_s)-(x_{i+1}+\ldots+x_s)$ is divisible by $p_{j_i+1}$, but both $x_i$, $y_i$ belong to $\mathrm{sum}(A_i)\subset \{0,1,\ldots,p_{j_i+1}-1\}$, a contradiction.
\end{proof}
\begin{proposition} \label{prop:distinct_irreducible_collections_have_different_spans}
    If two irreducible $\{d_i\}$-collections $A,B$ are distinct, then $\mathrm{sum}(A)\neq \mathrm{sum}(B)$.
\end{proposition}
\begin{proof}
Let $A=(a_0,\ldots,a_k)$, $B=(b_0,\ldots,b_l)$. Assume that $k\ne l$, say, $l>k$. Then $\mathrm{sum}(B)\geq p_l+\sum_{i=0}^{l-1} (d_i-1)p_i=2p_l-1$ while $\mathrm{sum}(A)\leq 2\sum_{i=0}^{k-1} (d_i-1)p_i = 2(p_{k+1})-1 \leq 2p_l-2 < \mathrm{sum}(B)$. If $l=k$, then choose the minimal index $i$ for which $a_i\ne b_i$. If $i=k$, then $\mathrm{sum}(A)\neq \mathrm{sum}(B)$ is clear. If $i<k$, then both $a_i,b_i$ belong to $\{d_i-1,\ldots,2d_i-2\}$, so $a_i$ and $b_i$ differ modulo $d_i$,
thus $\mathrm{sum}(A)$ and $\mathrm{sum}(B)$ differ modulo $p_i$.
\end{proof}

\begin{proposition}\label{prop:equal_spans_yield_equal_collections_for_normal}
If $A$ and $B$ are two normal $\{d_i\}$-collections, and $\mathrm{span}(A)=\mathrm{span}(B)$, then $A=B$. 
\end{proposition}
\begin{proof} Let $A=A_1\sqcup A_2\sqcup \ldots \sqcup A_s$ be a decomposition of $A$, $B=B_1\sqcup B_2\ldots \sqcup B_t$ that of $B$. The minimal non-element of $\mathrm{span}(A)$ equals $\mathrm{sum}(A_1)+1$ by Proposition \ref{prop:smallest_value_not_in_span} (if $s=1$, apply Lemma  \ref{lemma:condition_for_assuming_all_values} instead). Analogously for $B$, so we get $\mathrm{sum}(A_1)=\mathrm{sum}(B_1)$ and thus $A_1=B_1$ by Proposition \ref{prop:distinct_irreducible_collections_have_different_spans}. Now if $t$ or $s$ equals 1 we get from $\mathrm{sum}(A)=\mathrm{sum}(B)$ that $s=t=1$ and $A=B=A_1$. If both $t,s$ are greater than 1, denote by $p_j$ the maximal element of $A_1$, note that $\mathrm{sum}(A_1)<p_{j+1}$
again by Proposition \ref{prop:smallest_value_not_in_span}. Therefore, if we consider the map $\psi:x\to \lfloor x/p_{j+1}\rfloor$
on the set $\mathrm{span}(A)=\mathrm{span}(A_1)+\mathrm{span}(A_2)+\ldots+\mathrm{span}(A_s)$, its value
does not depend on the component in $\mathrm{span}(A_1)$, and we get 
$\psi(\mathrm{span}(A))=\frac{1}{p_{j+1}}(\mathrm{span}(A_2)+\ldots+\mathrm{span}(A_s))$. Use the same for $B$ to get
$\mathrm{span}(A_2)+\ldots+\mathrm{span}(A_s)=\mathrm{span}(B_2)+\ldots+\mathrm{span}(B_t)$ and proceed by induction
to conclude that $A_2=B_2$ and so on.
\end{proof}

\begin{theorem} \label{thm:spans_equal_iff_normalizations_equal}
    $\mathrm{span}(A)=\mathrm{span}(B)$ if and only if $\mathrm{norm}(A)=\mathrm{norm}(B)$
\end{theorem}
\begin{proof}
    Proposition \ref{prop:same_normalization_implies_same_span} and Proposition \ref{prop:equal_spans_yield_equal_collections_for_normal} 
    are the two directions.
\end{proof}

\begin{corollary} Start with arbitrary $\{d_i\}$-collection $A$ and apply elementary exchanges $e_i$ if $a_i\geq 2d_i-1$, in arbitrary order. Then this process terminates after finitely many steps, and the resulting normalized $\{d_i\}$-collection does not depend on the order.
\end{corollary}
\begin{proof} Note that we may get only finitely many
different $\{d_i\}$-collections, since the sum does not change. And there is a semiinvariant, for example, the sum of squares of elements,
which increases after every operation. Thus the process terminates. 
By Proposition \ref{prop:invariance_under_elementary_exchange}, these operations do not change the
span. Since normalized $\{d_i\}$-collection is uniquely determined by its span, the independence of the order follows.
\end{proof}

We conclude with the criterion used to determine if an elementary exchange breaks the desired span invariant.
\begin{proposition}
    $\mathrm{span}(e_i(A))=\mathrm{span}(A)$ if and only if $[\mathrm{norm}_i(A)]_i >2(d_i-1)$
\end{proposition}
\begin{proof}
    $\Rightarrow$: Suppose $[\mathrm{norm}_i(A)]_i \leq 2(d-1)$, in which case a run of Algorithm \ref{alg:proper_elem_exchange_seq} would not perform any exchanges at $i$, which means that $[\mathrm{norm}(A)]_i=A_i=[\mathrm{norm}_i(A)]_i > [\mathrm{norm}(e_i(A))]_i$, from which follows by one direction of Theorem \ref{thm:spans_equal_iff_normalizations_equal} that the spans differ.

    $\Leftarrow$: If $[\mathrm{norm}_i(A)]_i >2(d_i-1)$, then Algorithm \ref{alg:proper_elem_exchange_seq} which returns $\mathrm{norm}(A)$ necessarily performs at least one elementary exchange at $i$. Upon $e_i(A)$, we have performed the first elementary exchange there doing so on the places less significant than $i$. Thus, once that algorithm run on $e_i(A)$ reaches $i$, the state of the $\{d_i\}$-collection will be same as the state of $\{d_i\}$-collection right after the first elementary exchange on $i$ when the algorithm is run on $A$. Thus, $\mathrm{norm}(e_i(A))=\mathrm{norm}(A)$ in this case, from which follows by the other direction of Theorem \ref{thm:spans_equal_iff_normalizations_equal} that the spans are the same.
\end{proof}

\section{Acknowledgements}
I would like to thank Professor Fedor Petrov of Saint-Petersburg State University in the Department of Mathematics and Computer Science for previewing this paper and further verifying the correctness of its proofs. Moreover, it occurred to me to research this problem after I read parts of \cite{andrews_eriksson_integer_partitions} on the closely related topic of integer partitions.

\bibliography{multiset_sums_mixed_system}
\bibliographystyle{alpha}
\end{document}